\documentstyle[amstex,amscd,12pt]{article}
\title{The Hodge Numbers of the Moduli Spaces of Vector 
Bundles over a Riemann Surface} 
\author{Richard Earl and Frances Kirwan}
\begin{document}
\setcounter{section}{-1}
\setcounter{page}{1}
\newtheorem{prop}{PROPOSITION}
\newtheorem{lem}[prop]{LEMMA}
\newtheorem{cor}[prop]{COROLLARY}
\newtheorem{thm}[prop]{THEOREM}
\newtheorem{guess}{CONJECTURE}
\newtheorem{REM}[prop]{Remark}
\newenvironment{rem}{\begin{REM} 
\normalshape}{\end{REM}}
\newcommand{\Q}{{\bf Q}}
\newcommand{\C}{{\cal C}}
\newcommand{\G}{{\cal G}}
\newcommand{\mnd}{{\cal M}(n,d)}
\newcommand{\mond}{{\cal M}_{\Lambda}(n,d)}
\newcommand{\HS}{H^{*}}
\newcommand{\HG}{H^{*}_{\G}}
\newcommand{\ar}{a_{r}}
\newcommand{\brk}{b_{r}^{k}}
\newcommand{\fr}{f_{r}}
\newcommand{\Z}{{\bf Z}}
\newcommand{\T}{\mbox{\bf \normalshape t}}
\newcommand{\D}{\mbox{\normalshape d}}
\newcommand{\Res}{\mbox{\normalshape Res}}
\newcommand{\ch}{\mbox{\normalshape ch}}
\newcommand{\td}{\mbox{\normalshape td}}
\newcommand{\mtc}{M_{\T}(c)}\newcommand{\g}{\bar{g}}
\newcommand{\ng}{{\cal M}_{\Lambda}^{g}(2,1)}
\maketitle

\section{Introduction}

Let $\mnd$ denote the moduli space of stable holomorphic vector 
bundles of coprime rank $n$ and degree $d$ over a fixed Riemann 
surface $\Sigma$ of genus $g \geq 2$. Let $\Lambda$ be a fixed line bundle 
over $\Sigma$ of degree $d$ and let $\mond \subset \mnd$ denote the 
space consisting of those bundles with determinant $\Lambda$. The 
spaces $\mnd$ and $\mond$ are nonsingular 
complex projective varieties whose geometry has been much studied. In 
particular Harder, Narasimhan, Desale and Ramanan first described in 
1975 an inductive method to determine the Betti numbers of $\mnd$ using 
number theoretic methods and the Weil conjectures 
\cite{HN,DR}. Subsequently in 1982 Atiyah and Bott \cite{AB} obtained the 
same inductive method using gauge theory.\\   

\indent In this note we give a similar inductive method for 
determining the Hodge-Poincar\'{e} polynomials of $\mnd$ and $\mond$, that is 
\[HP(\mnd)(x,y) = \sum_{p \geq 0} \sum_{q \geq 0} h^{p,q}(\mnd) x^{p} y^{q},\] 
and 
\[HP(\mond)(x,y) = \sum_{p \geq 0} \sum_{q \geq 0} h^{p,q}(\mond) x^{p} y^{q},\] 
where $h^{p,q}$ denote the Hodge numbers (see Theorem \ref{thm} and 
Lemma \ref{lem} below). The $\chi(t)$-characteristic of $\mond$, that 
is $\chi(t) = HP(\mond)(t,-1),$ has an especially simple form (see 
Corollary \ref{chi}), while the $\chi(t)$-characteristic of $\mnd$ is 
identically zero.\\  

\indent In \cite{AB} Atiyah and Bott identify 
$\mnd$ with the quotient of ${\cal C}^{s}$, the infinite 
dimensional space of stable 
holomorphic structures on a fixed $C^{\infty}$ complex bundle $\cal E$ 
of rank $n$ and degree $d$ over $\Sigma$, by ${\cal G}_{c}$, the 
infinite dimensional group of 
smooth complex automorphisms. They introduce a stratification for 
${\cal C}$, the infinite dimensional affine space of all holomorphic 
structures on $\cal E$, which is equivariantly perfect with respect to 
the action of  ${\cal G}_{c}$ (or equivalently the gauge group $\cal G$) 
and which has ${\cal C}^{s}$ as an open stratum. The 
resulting Morse equalities give the equivariant cohomology of the 
stable stratum and thus the cohomology of $\mnd$, in terms of the 
classifying space of $\cal G$ and the equivariant cohomology of the 
unstable strata, which can be 
calculated inductively.\\   

\indent In \cite{K} the methods of Atiyah and Bott are adapted 
to finite dimensional quotients in the sense of Mumford's geometric 
invariant theory (GIT) \cite{M,N}, considering the linear action of a complex 
reductive group on a nonsingular complex projective variety where every 
semistable point is stable. This method generalises 
to give an inductive approach determining the Hodge numbers 
\cite[$\S$14]{K}.\\   

\indent To apply this to the case of $\mnd$ we recall from Newstead 
\cite{N} how this moduli space may be expressed as a finite dimensional 
quotient in the sense of GIT. It is shown in \cite{K2} that the 
resulting finite dimensional stratification coming from \cite{K} 
corresponds naturally to the 
stratification of Atiyah and Bott outside a subset whose codimension 
tends to infinity with $d$, and the 
equivariant cohomology of the corresponding strata agrees up to a 
degree tending to infinity with $d$. Since $\mnd$ depends on $d$ only 
through its remainder modulo $n$, we may then refine Atiyah and Bott's 
inductive formulas for the Betti numbers to give the Hodge numbers as well.\\   

\indent In the case $n=2$ our formula for the Hodge numbers of 
$\mond$ was recently proved by del Ba\~{n}o Rollin \cite[$\S$ 3]{BR} 
and had been independently discovered by Newstead \cite{N3}. The 
formula for the $\chi(t)$-characteristic of ${\cal M}_{\Lambda}(2,1)$ 
was first proved by Narasimhan and Ramanan \cite{NR} in 1975.\\   

\indent The layout of this paper is as follows. In Section 1 we review 
the arguments used in \cite{AB} and \cite{K2} to obtain inductive 
formulas for the Betti numbers of $\mond$. In Section 2 we adapt those 
arguments to Hodge numbers and prove our main result Theorem 
\ref{thm}. In Section 3 we explicitly calculate $HP(\mond)(x,y)$ for 
$n=2,3$. In Section 4 we discuss the $\chi(t)$-characterstic and 
explain how $\chi(t)$ contains information about intersection pairings 
involving certain generators of $\HS(\mond)$.  

\section{The Two Approaches} 

Let $\cal E$ be a fixed $C^{\infty}$ complex bundle of rank 
$n$ and degree $d$ over $\Sigma$ and let $\cal C$ denote the infinite 
dimensional affinespace of holomorphic structures on $\cal E$. For each holomorphic 
bundle $E \in {\cal C}$ there is a strictly ascending canonicalfiltration \cite[p.221]{HN} 
 \[0 = E_{0} \subset E_{1}  \subset \cdots \subset E_{P-1} \subset E_{P}= E\] 
such that the quotients $Q_{j} = E_{j}/E_{j-1}$ are semistable and such that 
\[\frac{\mbox{deg}(Q_{j})}{\mbox{rk}(Q_{j})} = \mu(Q_{j})>
 \mu(Q_{j+1}) = \frac{\mbox{deg}(Q_{j+1})}{\mbox{rk}(Q_{j+1})}.\] 
We then say that $E$ has type  
\begin{equation}\mu = ( \mu(Q_{1}), \ldots , 
\mu(Q_{P})) \in {\bf Q}^{n} \label{19}\end{equation} 
where $\mu(Q_{j})$ appears $\mbox{rk}(Q_{j})$ times. Atiyah and Bott 
defined a stratification $\{{\cal C}_{\mu}: \mu \in {\cal M} \}$ by 
setting ${\cal C}_{\mu} \subset {\cal C}$ to be the set of all 
holomorphic bundles $E$ over $\Sigma$ of type $\mu$. 
The set ${\cal C}^{ss}$ of semistable 
bundles is precisely the stratum of type 
\[\mu_{0} = \left( \frac{d}{n}, \cdots , \frac{d}{n} \right).\] 
\indent This stratification $\{{\cal C}_{\mu}: \mu \in {\cal M} \}$ is 
equivariantly perfect with respect to the action of the gauge group 
$\cal G$. Thus if we let 
\[P_{G}(X,t) = \sum_{j \geq 0} t^{j} \dim_{{\bf Q}}H^{j}_{G}(X,{\bf Q})\] 
denote the equivariant Poincar\'{e} polynomial of a space $X$ acted on 
by a group $G$, then we have equivariant Morse equalities 
\[P_{{\cal G}}({\cal C})(t) = P_{{\cal G}}({\cal C}^{ss})(t) + \sum_{\mu\neq \mu_{0}} 
t^{2d_{\mu}} P_{{\cal G}}({\cal C}_{\mu})(t)\] 
where $d_{\mu}$ is the complex codimension of ${\cal C}_{\mu}$ in $\cal C$ which 
is given by the formula \cite[7.16]{AB} 
\begin{equation}d_{\mu} = \sum_{1 \leq j < i 
\leq P} n_{i}d_{j} - n_{j}d_{i} +n_{i}n_{j}(g-1). \label{2}\end{equation} 
The ${\cal G}$-equivariant cohomology of the strata may be determined 
inductively via the isomorphisms \cite[7.12]{AB} 
\[H^{*}_{{\cal G}}({\cal C}_{\mu}) \cong \bigotimes_{1 \leq j \leq P} 
H^{*}_{{\cal G}(n_{j},d_{j})}({\cal C}(n_{j},d_{j})^{ss}).\]  
The $\G$-equivariant Poincar\'{e} polynomial of the affine space $\C$  
is given by \cite[Thm. 2.15]{AB}  
\begin{equation}P_{\G}(\C)(t) = \frac{\prod_{l=1}^{n} (1+t^{2l-1})^{2g}}{(1-t^{2n})
\prod_{l=1}^{n-1}(1-t^{2l})^{2}}. \label{9}\end{equation}  
Using this Atiyah and Bott obtain an inductive formula  
\begin{equation}P_{\G}(\C^{ss})(t) = P_{\G}(\C)(t) - \sum_{\mu\neq \mu_{0}} 
t^{2d_{\mu}} P_{\G(n_{j},d_{j})}(\C(n_{j},d_{j})^{ss})(t) \label{3}\end{equation}  
(where the sum is over all unstable types $\mu = (d_{1}/n_{1}, 
\ldots ,d_{P}/n_{P})$) 
for the $\G$-equivariant Betti numbers of  
$\C^{ss}$. They show that when $n$ and $d$ are coprime the quotient  
$\overline{\G}= \G/S^{1}$ of $\G$ by its central subgroup, consisting  
of multiplication by scalars in $S^{1}$, acts freely on $\C^{ss}$ and  
that \cite[9.3]{AB}  
\[P(\mnd)(t) = P_{\overline{\G}}(\C^{ss})(t) = (1-t^{2})P_{\G}(\C^{ss})(t).\]  
Finally they show that \cite[Prop. 9.7]{AB}  
\[P(\mnd)(t) = (1+t)^{2g} P(\mond)(t).\]  
\indent The GIT approach described in \cite{N} involves the action of  
a projective general linear group on a variety $R$ which naturally 
parametrises holomorphic  
bundles over $\Sigma$ as follows. For any stable bundle $E$ of rank  
$n$ and degree $d$ over $\Sigma$ of slope $d/n > 2g-1$ we 
have from \cite[5.2]{N} and the Riemann-Roch  
theorem that  
\begin{equation}E \mbox{ is generated by its 
sections} \label{cond1}\end{equation}  
and that  
\begin{equation}\dim H^{0}(\Sigma,E) = d+n(1-g) \qquad \mbox{and} 
\qquad H^{1}(\Sigma,E) =0. \label{cond2}\end{equation}  
Let 
\[p = d + n(1-g),\]  
and let $G(n,p)$ denote the Grassmannian of $n$-dimensional quotients  
of ${\bf C}^{p}$. There is then a natural one-to-one correspondence  
between holomorphic maps 
\[f:\Sigma \to G(n,p)\]  
and holomorphic quotients $E$ of $\Sigma \times {\bf C}^{p}$ of rank  
$n$, defined by setting the fibre $E_{x}$ to be  $f(x)$; more  
precisely $E = f^{*}Q$ where $Q$ is the tautological bundle of rank $n$  
over the quotient Grassmannian $G(n,p)$.\\  
\indent Let $\mbox{Hol}_{d}(\Sigma,G(n,p))$ denote the space of  
holomorphic maps $f:\Sigma \to G(n,p)$ such 
that the induced bundle $E= f^{*}Q$ has degree $d$. 
We then define $R\subset \mbox{Hol}_{d}(\Sigma,
G(n,p))$ to consist of those holomorphic maps 
such that $E=f^{*}Q$ satisfies $H^{1}(\Sigma,E) = 0$ and the map on  
sections \[{\bf C}^{p} \to H^{0}(\Sigma,E)\] is surjective. R is then a 
nonsingular quasi-projective variety 
\cite[p.265]{K2}.\\  
\indent We now recall some properties of the action of $PGL(p)$ 
on $R$and define a linearisation for this 
action. Assume that $d/n > 2g-1$. Then there is a quotient  
$\tilde{E}$ of the rank $p$ trivial complex 
bundle over $R \times\Sigma$ such that\\  
\indent (a) $\tilde{E}$ has the local universal property for families  
of bundles over $\Sigma$ of rank $n$ and degree $d$ which satisfy  
(\ref{cond1}) and (\ref{cond2}),\\  
\indent (b) for every $f \in R$ the restriction   
\[E^{f} = \tilde{E} \left|_{ \{ f \} \times \Sigma} \right. = f^{*}Q\]  
is the induced quotient of $\Sigma \times {\bf C}^{p}$,\\  
\indent (c) for $f,g \in R$ we have $E^{f} \cong E^{g}$ if and only if  
$f$ and $g$ lie in the same orbit of the action of $PGL(p)$ on $R$.\\  
\indent Now by \cite[pp. 141-3]{N} if $N$ is suitably large then $R$ can be  
embedded as a quasi-projective subvariety of $(G(n,p))^{N}$ by a map  
\[f \mapsto (f(x_{1}),...,f(x_{N}))\]  
where $x_{1},...,x_{N}$ are points of $\Sigma$. This provides the  
required linearisation for the action of $SL(p)$ on $R$. Moreover for  
sufficiently large $N$ and $d$ we have\\  
\indent (d) $f \in R^{ss}$ for this linearisation if and only if  
$E^{f}$ is semistable; further the closure of $R^{ss}$ in 
$(G(n,p))^{N}$ is in $R$.\\  
\indent When $n$ and $d$ are coprime integers, semistable bundles of  
rank $n$ and degree $d$ are stable 
and hence simple and thus $PGL(p)$ acts freely on $R^{ss}$. From GIT  
the quotient $R^{ss}/PGL(p)$ is a nonsingular complex projective 
variety which we have 
identified with $\mnd$. Thus the Poincar\'{e} polynomial of 
$\mnd$ is given by \cite[9.4]{K2}  
\begin{equation}P_{PGL(p)}(R^{ss})(t) = (1-t^{2})P_{GL(p)}
(R^{ss})(t) \label{5}\end{equation}  

\indent It is proved in \cite[Lemma 10.1]{K2} and \cite[Cor. 7.4]{K2}  that  
for any positive integer $k$ there is a natural isomorphism  
\begin{equation}H^{j}_{GL(p)}(R) \cong H^{j}_{\G}(\C) \label{14}\end{equation}  
for $j \leq k,$ provided that the degree $d$ is sufficiently large  
(depending on $k$ and $n$). We are always free to assume that the 
degree $d$ is arbitrarily large, because tensoring with a fixed line  
bundle over $\Sigma$ of degree $e$ induces an isomorphism  
\begin{equation}\mnd \to {\cal M}(n,d+ne) \label{1}.\end{equation}  
Note that tensoring with such a line bundle takes semistable bundles  
to semistable bundles and takes bundles of type  
\[\mu = (d_{1}/n_{1}, \ldots , d_{P}/n_{P})\]  
to bundles of type 
\[\mu + e =  (d_{1}/n_{1}+e, \ldots , d_{P}/n_{P}+e).\]  
Note also that by (\ref{9}) above, $P_{\G}(\C)(t)$ is independent of  $d$.\\   

\indent By \cite[$\S$ 11]{K2} there is a stratification 
$\{S_{\beta}:\beta \in {\cal B}\}$ of $R$ which is the intersection  
with $R$ of the stratification of the product of Grassmannians  
$(G(n,p))^{N}$ described in \cite[$\S$ 16]{K} and is closely related  
to the stratification $\{\C_{\mu}: \mu \in {\cal M}\}$ described  
above. In fact, to every possible type $\mu$ of a holomorphic bundle  
of rank $n$ and degree $d$ we can associate an index 
$\beta(\mu) \in{\cal B}$ \cite[11.1]{K2}. 
Moreover given any finite set $\cal U$ of types 
of bundles of rank $n$ and degree $d$, we can 
choose $N$ and $e$ sufficiently large that if 
$\mu \in {\cal U}+e$ then a point $f$ of $R$  
lies in the stratum indexed by $\beta(\mu)$ if and only if the bundle  
$E^{f}$ is of type $\mu$ \cite[Cor. 1.5]{K2}. We can also assume  
\cite[Lemma 11.6]{K2}, given any positive 
integer $k$ as above, that every stratum of $R$ not  
indexed by $\beta(\mu)$ for any $\mu \in {\cal U}+e$ has codimension  
greater than $k$ and that every stratum indexed by $\beta(\mu)$ for  
some $\mu \in {\cal U}+e$ has codimension the integer $d_{\mu}$  
defined at (\ref{2}) above, and equivariant Poincar\'{e} polynomial  
\[P_{GL(p)}(S_{\beta(\mu)})(t) = \prod_{j=1}^{P} P_{GL(p_{j})}(R(n_{j},d_{j})^{ss} )\]  
where $\mu = (d_{1}/n_{1},\ldots,d_{P}/n_{P})$ and $p_{j} = d_{j} +n_{j}(1-g)$ 
\cite[Lemma 12.3]{K2}. Finally by \cite[$\S$12]{K2} we can also assume that, 
outside a subset of codimension $k$, the stratification 
$\{S_{\beta}:\beta \in {\cal B}\}$ satisfies the 
conditions of \cite[9.5]{K} to be equivariantly 
perfect. Putting all this together, we obtain  
\begin{equation}P_{\G}(\C)(t) = P_{GL(p)}(R^{ss})(t) 
+ \sum_{\mu \neq \mu_{0}}t^{2d_{\mu}} 
\prod_{j=1}^{P} P_{GL(p_{j})} (R(n_{j},d_{j})^{ss})(t) 
+O(t^{k(n,d)}) \label{13}\end{equation}  
where $k(n,d) \to \infty$ as $d \to \infty$ and $n$ remains fixed 
\cite[$\S$13]{K2}. This gives us an inductive formula which enables  
us to calculate the $GL(p)$-equivariant Betti numbers of $R^{ss}$, and  
thus the Betti numbers of $\mnd$ when $n$ and $d$ are coprime, up to  
any arbitrarily large degree by using the isomorphism (\ref{1}). It  
is, of course, equivalent to the inductive formula (\ref{3}), but  
because its derivation involves stratifications of finite dimensional  
quasi-projective varieties instead of infinite dimensional spaces, we  
shall be able to apply \cite[$\S$14]{K} to generalise it to cover  
Hodge numbers as well as Betti numbers.  

\section{Hodge Numbers}  

In the last section we described inductive formulas from which the  
Betti numbers of the moduli spaces $\mnd$ and $\mond$ can be  
calculated. We now want to show that these formulas can be refined to 
enable us to calculate the Hodge numbers of $\mnd$ and $\mond$, as is  
done for GIT quotients of nonsingular projective varieties 
in \cite[$\S$ 14]{K}.  
As in \cite[$\S$ 14]{K} we have to make use of Deligne's extension of  
Hodge theory to complex quasi-projective varieties $Y$ which may be 
non-compact and singular \cite{D1,D2}. If Y is such a variety then its  
cohomology groups $H^{j}(Y;{\bf C})$ have two canonical filtrations,  
the weight filtration and the Hodge filtration, giving a mixed Hodge  
structure on each $H^{j}(Y;{\bf C})$. The Hodge numbers  
$h^{p,q}(H^{j}(Y))$ of $H^{j}(Y)$ are then defined to be the  
dimensions of appropriate quotients associated to these filtrations  
\cite[II 2.3.7]{D1} and satisfy  
\[\dim H^{j}(Y;{\bf C}) = \sum_{p,q} h^{p,q}(H^{j}(Y)).\]  
Moreover $h^{p,q}(H^{j}(Y))=0$ unless $p$ and $q$ lie between  
$\max\{0,j-m\}$ and $\min\{j,m\}$, where $m$ is the complex 
dimension of $Y$,  
and $p+q \leq j$ if $Y$ is projective whereas $p+q \geq j$ 
if $Y$ is nonsingular \cite[III 8.2.9]{D1}. 
If $Y$ is both projective and nonsingular then the Hodge 
numbers $h^{p,q}(H^{j}(Y))$ for $p+q=j$ are  
the classical Hodge numbers $h^{p,q}(Y)$.  Let us define 
the Hodge-Poincar\'{e} polynomial of a quasi-projective  
variety $Y$ as 
\[HP(Y)(x,y,t) = \sum_{p,q,j} x^{p} y^{q} t^{j} h^{p,q}(H^{j}(Y)).\]  
Of course when $Y$ is projective and nonsingular we lose no  
information by setting $t=1$ and defining  
\[HP(Y)(x,y) = \sum_{p,q} x^{p} y^{q} h^{p,q}(H^{p+q}(Y));\]  
we then have 
\begin{equation}HP(Y)(x,y,t) = HP(Y)(xt,yt). \label{label}\end{equation} 
We shall omit the variable $t$ when (\ref{label}) is satisfied. 
We can also define equivariant Hodge numbers $h^{p,q}(H^{j}_{G}(Y))$  
when $Y$ is acted on algebraically by a complex reductive group $G$,  
by identifying $H^{*}_{G}(Y)$ with $H^{*}(Y \times_{G} EG)$ where $EG  
\to BG$ is a universal classifying bundle for $G$. Although $EG$ and  
$BG$ are not finite dimensional manifolds, there are natural Hodge  
structures on their cohomology (which in the case of $EG$ is of course  
trivial) (see \cite[III $\S$9]{D1}). We can regard $EG$ and $BG$ as  
increasing unions of finite dimensional  
varieties $(EG)_{m}$ and $(BG)_{m}$ for $m \geq 1$ such that $G$ acts  
freely on $(EG)_{m}$ with $(BG)_{m} = (EG)_{m}/G$ and the inclusions  
of $(EG)_{m}$ and $(BG)_{m}$ induce isomorphisms of cohomology in  
degrees less than $m$ which preserve the Hodge structure. For example  
\[B{\bf C}^{*} = {\bf CP}_{\infty} = \bigcup_{m \geq 1} {\bf CP}_{m}\]  
has Hodge-Poincar\'{e} polynomial  
\begin{equation}HP(B{\bf C}^{*})(x,y) = \sum_{j \geq 0} x^{j} y^{j} 
=\frac{1}{1-xy}. \label{4} \end{equation}  
Similarly $Y \times_{G} EG$ is an increasing union of finite  
dimensional varieties whose mixed Hodge structures induce a natural  
mixed Hodge structure on the cohomology groups of $Y \times_{G} EG$.  
We want to generalise to Hodge numbers the arguments described in the  
last section to obtain an inductive formula giving the Betti numbers  
of $\mnd$ using the identification of $\mnd$ with the quotient  
$R^{ss}/GL(p)$. First recall that the equality (\ref{5}) came from  
isomorphisms (cf. \cite[$\S$ 9]{K2} and \cite[$\S$ 9]{AB})  
\begin{equation}H^{*}_{PGL(p)}(R^{ss}) \cong \HS(\mnd) \label{10}
\end{equation}  
and  
\begin{equation}H^{*}_{GL(p)}(R^{ss}) \cong  \HS(\mnd) \otimes 
\HS(B{\bf C}^{*}).\label{11}\end{equation}  
The first of these isomorphisms comes from the fact that the central  
subgroup ${\bf C}^{*}$ of $GL(p)$ acts trivially on $R^{ss}$ and the  
induced action of the quotient $PGL(p) = GL(p)/{\bf C}^{*}$ is  
free. The natural map\[R^{ss} \times_{PGL(p)} EPGL(p) \to 
R^{ss}/PGL(p) \cong \mnd\]  
is a fibration with contractible fibre $EPGL(p)$ and so induces the  
isomorphism (\ref{10}). As morphisms of nonsingular quasi-projective  
varieties induce maps on cohomology which are strictly compatible with  
both the Hodge filtration and the weight filtration \cite[II3.2.11.1]{D1}, the 
isomorphism (\ref{10}) respects the Hodgestructures.  
The natural map 
\[R^{ss} \times_{GL(p)} EGL(p) \to R^{ss}/GL(p) \cong \mnd\]  
is a fibration with fibre $EGL(p)/{\bf C}^{*}$ which can be  
identified with $B{\bf C}^{*}$ since $EGL(p)$ is a contractible space  
on which ${\bf C}^{*}$ acts freely. This fibration is cohomologically  
trivial \cite[9.3]{AB} and so induces an isomorphism between 
$H^{*}_{GL(p)}(R^{ss})$ 
and $\HS(\mnd) \otimes\HS(B{\bf C}^{*})$, which is an 
isomorphism of Hodge structures  
(cf. \cite[III 8.1]{D1} and \cite[III 8.2.10]{D1}). Thus by (\ref{4})we have 
\begin{equation}HP_{GL(p)}(R^{ss})(x,y) = \frac{1}{1-xy} HP(\mnd)(x,y). 
\label{12}\end{equation}  
Next let us consider, from the point of view of Hodge numbers, the  
isomorphism (\ref{14}) between $\HS_{GL(p)}(R)$ and $H^{*}_{\G}(\C)$  
up to some degree tending to infinity with $d$. Since $\C$ is  
contractible $H^{*}_{\G}(\C)$ is canonically isomorphic to $\HS(B\G)$  
where $B\G$ is the classifying space for the gauge group $\G$. By  
\cite[Prop. 2.4]{AB} $B\G$ can be identified with a component  
$\mbox{Map}_{d}(\Sigma,BU(n))$ of the space of continuous maps  
$\mbox{Map}(\Sigma,BU(n))$ from $\Sigma$ to $BU(n)$. By  
\cite[Thm. 2.15]{AB} and \cite[Prop. 2.20]{AB} $\HS(B\G)$ is freely  
generated as a polynomial algebra tensored with an exterior algebra by  
generators  
\[a_{r} \in H^{2r}(B\G), \quad b_{r}^{j} \in H^{2r-1}(B\G),\]  
for $1 \leq r \leq n$ and $1 \leq j \leq 2g$, and by  
\[f_{r} \in H^{2r-2}(B\G)\]  
for $2 \leq r \leq n,$ where 
\[a_{r} \otimes 1 + \sum_{j=1}^{2g} b_{r}^{j} \otimes 
\alpha_{j} + f_{r}\otimes \omega\]  
is the $r$th Chern class of the universal 
bundle $\cal{V}$ over $B\G\times \Sigma$ and 
$\{1\},\{\alpha_{j}: 1 \leq j \leq 2g\}$ and $\{\omega\}$ 
are the standard bases for 
$H^{0}(\Sigma),H^{1}(\Sigma)$ and $H^{2}(\Sigma)$. 
The universal bundle ${\cal V}$ is 
the pullback to $B\G \times \Sigma$ of the universal bundle over the  
infinite Grassmannian $BU(n)$ via the evaluation map  
\[B\G \times \Sigma = \mbox{Map}_{d}(\Sigma,BU(n)) 
\times \Sigma \to BU(n).\]  
We can take the total space $EGL(p)$ for the classifying bundle for  
$GL(p)$ to be the space of all surjective linear maps from  
${\bf C}^{\infty}$ to ${\bf C}^{p}$ with the obvious action of  
$GL(p)$ (see \cite[$\S$ 7]{AB} or \cite[$\S$ 3]{VLP} for further  
details). If we similarly identify $BU(n)$ with the Grassmannian  
$G(n,\infty)$ of $n$-dimensional quotients (or equivalently subspaces  
of codimension $n$) of ${\bf C}^{\infty}$, then we obtain a map  
\[\theta: R \times_{GL_{p}} EGL(p) \to  \mbox{Map}_{d}(\Sigma,BU(n)) 
=  \mbox{Map}_{d}(\Sigma,G(n,\infty))\]  
defined as follows. Given an element $f: \Sigma \to G(n,p)$ of $R$ and  
an element $e:{\bf C}^{\infty} \to {\bf C}^{p}$ of $EGL(p)$ then  
$\theta(f,e)$ is the map from $\Sigma$ to $G(n,\infty)$ which 
sends $x\in \Sigma$ 
to the kernel of the composition of the surjection $e:{\bf C}^{\infty} 
\to {\bf C}^{p}$ with the projection of  
${\bf C}^{p}$ onto its quotient $f(x) \in G(n,p)$. By \cite[Lemma10.1]{K2} 
and \cite[Cor. 7.4]{K2} this map 
$\theta$ induces an isomorphism  
\[H^{*}_{\G}(\C) = \HS(B\G) \to H^{*}_{GL(p)}(R)\]  
up to some degree tending to infinity with $d$, as at (\ref{14})  
above. The pullback via $\theta$ of the universal bundle $\cal{V}$  
over $B\G \times \Sigma$ corresponds to the pullback via the  
evaluation map  
\[R \times \Sigma \subseteq \mbox{Hol}_{d}(\Sigma,G(n,p)) 
\times \Sigma\to G(n,p)\]  
of the tautological bundle of rank $n$ over the Grassmannian $G(n,p)$  
of $n$-dimensional quotients of ${\bf C}^{p}$. It is therefore a  
holomorphic bundle over $R \times \Sigma$ and its Chern classes are of  
type $(p,p)$ \cite[p.417]{GH}. Now the generators $a_{r},b_{r}^{j}$  
and $f_{r}$ of $\HS(B\G)$ were defined by the equation  
\[c_{r}(\cal{V}) =a_{r} \otimes 1 + \sum_{j=1}^{2g} b_{r}^{j} \otimes\alpha_{j} 
+ f_{r} \otimes \omega\]  
for $c_{r}(\cal{V}) \in \HS(B\G) \otimes \HS(\Sigma)$, where 
$\omega\in H^{2}(\Sigma)$ 
has Hodge type (1,1) and we may choose the basis  
$\{\alpha_{j}:1 \leq j \leq 2g\}$ for $H^{1}(\Sigma)$ such that  
$\alpha_{j}$ has type (1,0) if $1 \leq j \leq g$ and Hodge type (0,1)  
if $g+1 \leq j \leq 2g$. Since $\theta^{*}(c_{r}({\cal V})) = 
c_{r}(\theta^{*}{\cal V})$ has Hodge type $(r,r)$ it follows that  
under the isomorphism (\ref{14}) the generators $a_{r}$ and $f_{r}$  
for $\HG(\C)= \HS(\C)$ correspond to elements of Hodge type $(r,r)$  
and $(r-1,r-1)$ respectively, whereas the generators $b_{r}^{j}$  
correspond to elements of Hodge type $(r-1,r)$ if $1 \leq j \leq g$  
and $(r,r-1)$ if $g+1 \leq j \leq 2g$. Thus  
\begin{equation}HP_{GL(p)}(R)(x,y) = \frac{\prod_{l=1}^{n} (1+ x^{l}y^{l-1})^{g} 
(1+x^{l-1}y^{l})^{g}}{(1-x^{n}y^{n})\prod_{l=1}^{n-1} 
(1-x^{l}y^{l})^{2}}+ O((xy)^{k(n,d)}) \label{17}\end{equation} 
where $k(n,d)$ tends to infinity with $d$.  
We can now prove  

\begin{thm} \label{thm}
The Hodge-Poincar\'{e} polynomial of the moduli space $\mnd$, when $n$ 
and $d$ are coprime, satisfies  
\begin{equation}HP(\mnd)(x,y) = (1-x y) F_{n,d}(x,y) \label{20}\end{equation}  
where $F_{n,d}(x,y) = HP_{GL(p)}(R^{ss})(x,y)$ is given by the  
inductive formula 
\begin{equation}\label{21} F_{n,d}(x,y) = \frac{\prod_{l=1}^{n} 
(1+ x^{l}y^{l-1})^{g} (1+x^{l-1}y^{l})^{g}}{(1-x^{n}y^{n})
\prod_{l=1}^{n-1} (1-x^{l}y^{l})^{2}} - \sum_{\mu \neq \mu_{0}} 
(xy)^{d_{\mu}} \prod_{1 \leq j \leq P}F_{n_{j},d_{j}}(x,y).
\end{equation}
Here the sum is over all types $\mu= (d_{1}/n_{1},\ldots ,d_{P}/n_{P})$  
of bundles of rank $n$ and degree $d$ as described at (\ref{19}), and  
$d_{\mu}$ is given by (\ref{2}).\end{thm}

\begin{rem} The inductive formula (\ref{21}) is to be used even  
when $n$ and $d$ are not coprime. It is for the identity  
(\ref{20}) that $n$ and $d$ need be coprime.\end{rem}  

{\bf PROOF:} We can use equation (\ref{21}) to define $F_{n,d}(x,y)$  
by induction on $n$, because when $n=1$ there are 
no types $\mu \neq\mu_{0}$. 
Then by (\ref{12}) above it suffices 
to prove that when $n$and $d$ are coprime 
\[HP_{GL(d+n(1-g))}(R(n,d)^{ss})(x,y) = F_{n,d}(x,y)\]  
where we have written $R(n,d)$ for $R$ and $d+n(1-g)$ for $p$ to make  
explicit the dependence on $n$ and $d$.\\  
\indent First note that 
\[\frac{\prod_{l=1}^{n} (1+ x^{l}y^{l-1})^{g} (1+
x^{l-1}y^{l})^{g}}{(1-x^{n}y^{n})\prod_{l=1}^{n-1} (1-x^{l}y^{l})^{2}}\]  
is independent of $d$ and that  
\[d_{\mu} = \sum_{1 \leq j<i \leq P} n_{i}n_{j}\left(\frac{d_{j}}{n_{j}} 
- \frac{d_{i}}{n_{i}} + g - 1 \right)\]  
is unchanged if $\mu = (d_{1}/n_{1},\cdots,d_{P}/n_{P})$ is replaced  
by $\mu + e = (d_{1}/n_{1} + e, \cdots, 
d_{P}/n_{P} + e)$ for any $e\in {\bf Z}$. 
Since $\mu+e$ runs over all types of holomorphic  
bundles of rank $n$ and degree $d+ne$ as $\mu$ runs over types of  
holomorphic bundles of rank $n$ and degree $d$ (see $\S$1), it  
follows by induction on $n$ that $F_{n,d}(x,y)$ is  
unchanged if $d$ is replaced by $d+ne$ for any $e \in {\bf Z}$.  
Now recall from $\S$1 that there is a stratification $\{S_{\beta}:
\beta \in {\cal B}\}$ of $R$ which is the intersection  
with $R$ of the stratification of the product of Grassmannians  
$(G(n,p))^{N}$ described in \cite[$\S$16]{K} and is  
$GL(p)$-equivariantly perfect outside a subset of dimension $k(n,d)$  
where $k(n,d) \to \infty$ as $d \to \infty$. By \cite[$\S$ 14]{K} the  
equivariant Morse equalities associated to this stratification can be  
refined to give the following relation between equivariant  
Hodge-Poincar\'{e} polynomials  
\begin{equation}HP_{GL(p)}(R)(x,y) = 
\sum_{\beta \in {\cal B}} (xy)^{d_{\beta}} 
HP_{GL(p)}(S_{\beta})(x,y) + O((xy)^{k(n,d)}) \label{22}\end{equation}  
where $d_{\beta}$ is the complex codimension of $S_{\beta}$ in $R$  
(see \cite[14.5]{K}). We noted in $\S$1 that to every possible type  
$\mu$ of a holomorphic bundle of rank $n$ and degree $d$ we can  
associate an index $\beta(\mu)$ such that every stratum of $R$ of  
codimension less than $k(n,d)$ is indexed by $\beta(\mu)$ for some  
type $\mu = (d_{1}/n_{1}, \ldots , d_{P}/n_{P})$. Moreover  
$S_{\beta(\mu)}$ then has codimension 
\begin{equation}d_{\beta_{\mu}} = d_{\mu} \label{23}\end{equation}  
where $d_{\mu}$ is defined by (\ref{2}), and equivariant cohomology  
\[\HS_{GL(p)}(S_{\beta(\mu)}) \cong \bigotimes_{j=1}^{P}
\HS_{GL(p_{j})} (R(n_{j},d_{j})^{ss})\]  
where $p_{j} = d_{j} + n_{j}(1-g)$. This isomorphism of equivariant  
cohomology is induced by inclusions 
of subvarieties and subgroups (see 
the proof of \cite[Lemma 12.3]{K2}) and therefore respects Hodge  
structures. Thus 
\begin{equation}HP_{GL(p)}(S_{\beta})(x,y) = \prod_{j=1}^{P} 
HP_{GL(p_{j})}(R(n_{j},d_{j})^{ss})(x,y). \label{24}\end{equation}  
Since $R^{ss}$ is the stratum indexed by $\beta(\mu_{0})$ where  
$\mu_{0} = (d/n,\ldots, d/n)$, combining (\ref{17}), (\ref{22}), 
(\ref{23}) and (\ref{24}) yields  
\begin{eqnarray*}HP_{GL(p)}
(R^{ss})(x,y) = \frac{\prod_{l=1}^{n} (1+ x^{l}y^{l-1})^{g} 
(1+x^{l-1}y^{l})^{g}}{(1-x^{n}y^{n})\prod_{l=1}^{n-1} 
(1-x^{l}y^{l})^{2}}\\- \sum_{\mu \neq \mu_{0}} 
(xy)^{d_{\mu}} \prod_{1 \leq j \leq P}HP_{GL(
d_{j}+n_{j}(1-g))}(R(n_{j},d_{j})^{ss})(x,y) + O((xy)^{k(n,d)}).
\end{eqnarray*}  
where $k(n,d) \to \infty$ as $d \to \infty$. By induction on $n$ it  
follows that, given any $k>0$, if $e$ is sufficiently large then  
\begin{eqnarray*}HP_{GL(d+ne +n(1-g))}(R(n,d+ne)^{ss})(x,y) & 
= & F_{n,d+ne}(x,y) +O((xy)^{k})\\& = & F_{n,d}(x,y) + O((xy)^{k})
\end{eqnarray*}  
But by (\ref{1}) and (\ref{12}), when $n$ and $d$ are coprime,  
$HP_{GL(d+ne+n(1-g))}(R(n,d+ne)^{ss})(x,y)$ is independent of $e$ and  
hence we must have\[HP_{GL(d+n(1-g))}(R(n,d)^{ss})(x,y) = F_{n,d}(x,y)\]  
as required. 
\indent $\Box$ \\  
[\baselineskip]  
\indent This theorem enables us 
to calculate the Hodge numbers of $\mnd$. It  
also gives us the Hodge numbers of the moduli space $\mond$ of bundles  
with fixed determinant $\Lambda$ via the following lemma  

\begin{lem}\label{lem}  
The Hodge-Poincar\'{e} polynomials of $\mnd$ and $\mond$ satisfy  
\[HP(\mnd)(x,y) = (1+x)^{g} (1+y)^{g} HP(\mond)(x,y).\]\end{lem}

{\bf PROOF:} Atiyah and Bott show \cite[Prop. 9.7]{AB} that the  
determinant map\[\mnd \to {\cal M}(1,d) = \mbox{Jac}_{d}\]  
from $\mnd$ to the Jacobian $\mbox{Jac}_{d}$ with fibre $\mond$ is a  
cohomologically trivial fibration and so induces an 
isomorphism of Hodge  
structures (cf. \cite[III 8.1]{D1} and \cite[III 8.2.10]{D1})  
\[\HS(\mnd) \cong \HS(\mond) \otimes \HS(\mbox{Jac}_{d}).\]  
The Hodge-Poincar\'{e} polynomial of the Jacobian is  
\[HP(\mbox{Jac}_{d})(x,y) = (1+x)^{g} (1+y)^{g},\]  
(for example by Theorem \ref{thm} above with $n=1$) 
and so the result follows. 
\indent $\Box$

\begin{rem} When $d$ is odd then a related formula 
for the dimensions of the Hodge  
cohomology groups\[H^{p,p}({\cal M}_{\Lambda}(2,1)) 
\cap  H^{2p}({\cal M}_{\Lambda}(2,1);\Q)\]  
was found by Balaji, King and Newstead \cite{BKN}. 
If, for any nonsingular projective variety $X$, we define  
\[P_{H}(X)(t) = \sum_{p \geq 0} t^p ( \dim_{\Q} H^{p,p}(X) 
\cap H^{2p}(X;\Q) )\]  
then they show that  
\begin{equation} \label{bkn}P_{H}({\cal M}_{\Lambda}(2,1))(t) = 
\frac{P_{H}(J)(t^3) - t^gP_{H}(J)(t)}{(1-t)(1-t^2)}\end{equation}  
where $J$ is the Jacobian of $\Sigma$. The proof of  
Theorem \ref{thm} and Lemma \ref{lem} above could be modified to give  
an alternative proof of this result and in principle to give an  
inductive formula for $P_{H}(\mond)(t)$ when $n>2$, although the  
latter would be more complicated and would involve looking at the  
Hodge cohomology of the products of moduli spaces.\\  

\indent In \cite{BKN} the main object of study was the algebraic  
cohomology of ${\cal M}_{\Lambda}(2,1)$. By proving the analogous  
formula to (\ref{bkn}) for the algebraic cohomology it was shown in  
\cite{BKN} that the Hodge conjecture for ${\cal M}_{\Lambda}(2,1)$ is  
valid provided that it holds for the Jacobian $J$, which is true for  
generic curves $\Sigma$. In \cite{BN} the Hodge conjecture for generic  
$\Sigma$ was proved directly for $\mond$ with $n \geq 2$ and more  
generally for moduli spaces of parabolic bundles, so for such $\Sigma$  
the algebraic cohomology groups of $\mond$ would be given by  
$P_{H}(({\cal M}_{\Lambda}(2,1))(t)$. \end{rem}
       
\section{Calculations}  
In this section we use the inductive formula of Theorem \ref{thm}, together  
with Lemma \ref{lem}, to give explicit formulas for the Hodge numbers  
of the moduli space $\mond$.  
\begin{cor}  The Hodge-Poincar\'{e} 
polynomials of  $\mond$ for the following low  
values of $n$ are:\\  

\indent (a) when $n=2$ and $d$ is odd,  
\[\frac{(1+x^{2}y)^{g} (1+x y^{2})^{g} - x^{g} 
y^{g} (1+x)^{g}(1+y)^{g}}{(1-x y)(1-x^{2}y^{2})},\]  

\indent (b) when $n=3$ and $d \equiv1,2$ (mod 3),\[  
\frac{1}{(1-xy)(1-x^{2}y^{2})^{2}(1-x^{3}y^{3})}\left( (1+x^{2}y^{3})^{g} 
(1+x^{3}y^{2})^{g} (1+xy^{2})^{g}(1+x^{2}y)^{g} \right. \]  
 \[- x^{2g-1}y^{2g-1} (1+x y)^{2} (1+x)^{g} (1+y)^{g} 
(1+x y^{2})^{g}(1+x^{2} y)^{g} \] 
\[ \left. +x^{3 g -1} y^{3 g -1} (1+x y + x^{2} y^{2}) 
(1+x)^{2g} (1+y)^{2g}\right).\]
\end{cor}

{\bf PROOF:} (a) For $n=1$ Theorem \ref{thm} shows that for any $d$,  
\[F_{1,d}(x,y) = \frac{(1+x)^{g}(1+y)^{g}}{1-xy}.\]  
When $n=2$ and $d=1$ then the unstable types $\mu\neq \mu_{0}$ 
may each be written as 
$\mu(r) = (r+1,-r)$ for $r \geq0$ and from (\ref{2}) we note that 
$d_{\mu(r)} = 2 r + g$. For $n=2$ 
and $d=1$equation (\ref{21}) states that $F_{2,1}(x,y)$ equals 
\[\frac{(1+ x)^{g} (1+y)^{g} (1+ x^{2}y)^{g} (1+x y^{2})^{g}}{(1-x^{2}
y^{2}) (1-xy)^{2}} - \sum_{r=0}^{\infty} 
(xy)^{2r+g} \frac{(1+x)^{2g} (1+y)^{2g}}{(1-xy)^{2}},\]  
which in turn equals  
\begin{equation}\frac{(1+ x)^{g} (1+y)^{g} (1+ x^{2}y)^{g} 
(1+x y^{2})^{g} - (xy)^{g} (1+x)^{2g} (1+y)^{2g} } 
{(1-x^{2}y^{2})(1-xy)^{2}}. \label{2,1}\end{equation}  
Applying Lemma \ref{lem} and (\ref{20}) together with 
the isomorphism\ref{1} yields (a) above.\\  
[\baselineskip]  
\indent (b) For the case when $n=2$ and $d=0,$ we note that the types  
$\mu \neq \mu_{0}$ are of the form $\mu(r) = (r,-r)$ for $r \geq 1$  
with $d_{\mu}= 2r+g-1$. Arguing as in (a) above we see that  
$F_{2,0}(x,y)$ equals 
\begin{equation}\frac{(1+ x)^{g} (1+y)^{g} (1+ x^{2}y)^{g} (1+x 
y^{2})^{g+1}
 - (xy)^{g} (1+x)^{2g} (1+y)^{2g} }{(1-x^{2}y^{2})(1-xy)^{2}}. 
\label{2,0}\end{equation}  
Now let $n=3$. Notice that the map $E \mapsto E^{*} 
\otimes L$, where $L$ is a degree one line 
bundle over $\Sigma$, induces an  
isomorphism between the Hodge structures 
of ${\cal M}(3,1)$ and ${\cal M}(3,2)$ so that we may 
set $d=1$ without loss of generality. For  
$n=3$ the types $\mu \neq \mu_{0}$ fall into three types which we shall 
refer to as (1,1,1), (2,1) and (1,2).\\  

\indent (i) The (1,1,1)-types are  
\[\mu =(d_{1},d_{2},d_{3}) \mbox{ where } 
d_{1}>d_{2}>d_{3} \mbox{ and }d_{1}+d_{2}+d_{3} = 1,\]  
and $d_{\mu}$ equals $2(d_{1} - d_{3}) + 3(g-1).$ Let $r=d_{1}-d_{2}$ and  
$s=d_{2}-d_{3}$. Then \[d_{1} = \frac{1}{3} (2r+s+1), 
\indent d_{2} = \frac{1}{3}(-r+s+1)\]  
and so we must have $r-s \equiv 1 \mbox{ mod } 3.$ So the possible  
values of $r$ and $s$ are given by  
\[(r,s) = \quad (3k+1,3l+3) \quad \mbox{ or } \quad (3k+2,3l+1) 
\quad\mbox{ or } \quad (3k+3,3l+2)\]  
for $k \geq 0,l \geq 0.$ Thus  
$\sum_{(1,1,1)} (xy)^{d_{\mu}} \prod_{i=1}^{3} F_{1,d_{i}}(x,y)$ equals  
\[\frac{(1+x)^{3g}(1+y)^{3g}}{(1-xy)^{3}} 
\sum_{r-s \equiv 1 (3)}^{\infty} (xy)^{2r+2s+3g-3}\]  
\[= \frac{(xy)^{3g-3} (1+x)^{3g}(1+y)^{3g}}{(1-xy)^{3}} 
\sum_{k=0}^{\infty}\sum_{l=0}^{\infty} [(xy)^{2(3k+3l+4)} 
+ (xy)^{2(3k+3l+3)} +(xy)^{2(3k+3l+5)}]\]  
which, with a little simplifying, becomes  
\begin{equation}\frac{(xy)^{3g+3} (1+x)^{3g}(1+y)^{3g}}{(1-
xy)^{3} (1-x^2 y^2) (1-x^6y^6)}. \label{(1,1,1)} \end{equation} 
\indent (ii) The (2,1)-types are $\mu(r) = (r/2,r/2,1-r)$ for 
$r \geq1$ and $d_{\mu(r)}$ equals $3r+2g-4$. Thus the 
contribution of the (2,1)-types to (\ref{21}) equals  
\[\frac{(1+x)^{g}(1+y)^{g}}{1-xy} \sum_{r=1}^{\infty} 
(xy)^{3r + 2 g -4} F_{2,r}(x,y).\]  
As $F_{2,0}= F_{2,2r}$ and $F_{2,1} = F_{2,2r-1}$ for $r \geq 1$ then 
this in turn equals  
\[\frac{(1+x)^{g}(1+y)^{g}}{1-xy} \left( F_{2,0}(x,y) \sum_{r=1}^{\infty} 
(xy)^{6r+2g-4}  
+ F_{2,1}(x,y) \sum_{r=1}^{\infty}  (xy)^{6r+2g-7} \right)\]  
\begin{equation}= \frac{(1+x)^{g}(1+y)^{g}}{(1-xy)(1-x^6 y^6)} 
[ (xy)^{2g+2}F_{2,0}(x,y) 
+ (xy)^{2g-1} F_{2,1}(x,y)]. \label{(2,1)}\end{equation}  
\indent (iii) Similarly the (1,2)-types are given by 
$\mu(r) =(r,(1-r)/2,(1-r)/2)$ for $r \geq 1$ and 
$d_{\mu(r)}$ equals$3r+2g-3$. Arguing as in 
part (ii) we see that the contribution of the  
(1,2)-types to (\ref{21}) equals  
\begin{equation}\frac{(1+x)^{g}(1+y)^{g}}{(1-xy)
(1-x^6 y^6)} [ (xy)^{2g}F_{2,0}(x,y) + (xy)^{2g+3} 
F_{2,1}(x,y)]. \label{(1,2)}\end{equation}  
\indent Substituting the terms (\ref{2,1}), 
(\ref{2,0}),(\ref{(1,1,1)}), (\ref{(2,1)}) and (\ref{(1,2)}) into (\ref{21}) we  
obtain the final expression (b). 
\indent $\Box$
    
\begin{cor}\label{chi}  
The $\chi(t)$-characteristic of 
$\mond$, that is $HP(\mond)(t,-1)$, equals  
\begin{equation}\left( \prod_{r=1}^{n-1} (1 - (-t)^{r}) 
(1-(-t)^{r+1})\right)^{g-1}. \label{chi2}\end{equation}  
The $\chi(t)$-characteristic of $\mnd$ is identically zero.\end{cor}  
{\bf PROOF:} Note from Lemma \ref{lem} and (\ref{20}) that 
\[F_{n,d}(x,y) = \frac{(1+x)^g (1+y)^g HP(\mond)(x,y)}{1-xy}.\]  
Substituting this into the recurrence relation (\ref{21}) we obtain  
\[HP(\mond)(x,y) = \frac{\prod_{l=2}^{n} (1+ x^{l}y^{l-1})^{g} 
(1+x^{l-1}y^{l})^{g}}{\prod_{l=2}^{n}(1-x^{l}y^{l}) 
\prod_{l=1}^{n-1}(1-x^l y^l)}\]     
\[- \sum_{\mu \neq \mu_{0}} (xy)^{d_{\mu}} 
\frac{(1+x)^{(P-1)g}(1+y)^{(P-1)g}}{(1-xy)^{(P-1)g}} 
\prod_{1 \leq j \leq P} HP({\cal M}_{\Lambda}(n_{j},d_{j}))(x,y).\]     
As $P \geq 2$ for each $\mu \neq \mu_{0}$, then setting $x=t$ 
and $y=-1$ yields the required 
result for $\mond$. That the $\chi(t)$-characteristic of $\mnd$ 
is identically zero follows from Lemma 
\ref{lem}. \indent $\Box$       

\begin{cor}The Euler characteristic and signature of $\mond$ are zero.
\end{cor}  

{\bf PROOF:} These are respectively $\chi(-1)$ and $\chi(1)$.    

\section{The $\chi(t)$-characteristic}   

Let $N = \dim_{\bold C}(\mond) = (n^{2}-1)(g-1)$ and let $c_{r}$  
denote the $r$th Chern class of $\mond$ for $0 \leq r \leq N$. We also  
introduce here the notation $\g = g-1$. By the  
Riemann-Roch Theorem the $\chi(t)$-characteristic may be 
expressed as an evaluation on the fundamental class $[\mond]$ of  
$\mond$. Let $T$ denote the tangent bundle of $\mond$. 
We then have\cite[p. 142]{S} that 
\begin{equation}\chi(t) = (-1)^{N} \left( \sum_{p=0}^{N} t^{p} 
\ch (\bigwedge{}\!\!^{N-p}\,T^{*}) \, 
\td T \right) [ \mond ] . \label{chi3}\end{equation}  
Thus the $\chi(t)$-characteristic which we determined 
in Corollary\ref{chi} contains information, 
though limited, on the Chern numbers of 
$\mond$. From (\ref{chi2}) we know that  
\begin{equation}\chi(-1) = c_{N}[\mond] = 0. 
\label{ultimate}\end{equation}  
More generally, for $2 \leq k \leq N$, we may write  
$\chi^{(k)}(-1)$ as a linear combination of Chern numbers each of  
which includes at least one $c_{i}$ with $i > N - 2[k/2]$. (See  
\cite[$\S$ 3]{S} for details.) For example, \cite[3.7]{S}  
\begin{equation}\chi''(-1) = \frac{1}{24}(2 c_{1} c_{N-1} + 
N (3 N - 5) c_{N})[\mond],\label{penultimate}\end{equation}  
although for large values of $k$ these combinations of Chern numbers  
rapidly become more complicated. From (\ref{ultimate}) we see that  
\[c_{N}=0.\]  
By (\ref{penultimate}) it follows that $c_{1}c_{N-1}=0$. In fact as  
$c_{1} = 2 n f_{2}$ and since $f_{2}$ generates $H^{2}(\mond)$ then  
\[c_{N-1}=0\]  
by Poincar\'{e} duality. Unfortunately it is not obvious whether we  
may deduce any further vanishing of the Chern classes from these  
relations. It was proved by Gieseker \cite{G} that $c_{r}=0$ when 
$n=2$ and $r>2g-2$ 
and it has been conjectured \cite{EK} for $n > 2$  that $c_{r}=0$ 
when $r > n(n-1)(g-1)$.\\ 

\indent We may rewrite the pairing (\ref{chi3}) in terms of pairings  
involving the generators $a_{r}, b_{r}^{s}$ and $f_{r}$ as follows. If  
$\gamma_{1}, \cdots, \gamma_{N}$ are the Chern roots of $T$ then 
(see\cite[p.430]{H}) 
\[\sum_{p=0}^{N} t^{p} \ch (\bigwedge{}\!\!^{N-p}\, T^{*}) = 
\prod_{p=1}^{N} (t+ e^{-\gamma_{p}}).\]  
More generally, if 
\[\ch T = \sum_{p=1}^{N} e^{\gamma_{p}} = \sum_{i=1}^{K} 
\mu_{i}e^{\epsilon_{i}},\]  
where $\mu_{1}, \cdots , \mu_{K}$ and $\epsilon_{1},\cdots,
\epsilon_{K}$ are formal 
cohomology classes of degrees zero and two respectively, then 
\begin{equation}\sum_{p=0}^{N} t^{p} \ch (\bigwedge{}\!\!^{N-p}\, T^{*}) 
= (1+t)^{N}\prod_{i=1}^{K} \left( 1 + \frac{e^{-\epsilon_{i}} - 1}{1+t}
\right)^{\mu_{i}}, \label{newarr}\end{equation} 
where $(1+y)^{\mu}$ is to be interpreted as the formal power series 
\[\exp(\mu \log (1+y)) = 1 + \mu y + \frac{\mu (\mu-1)}{2!} y^2 + 
\frac{\mu (\mu-1)(\mu-2)}{3!} y^3 + \cdots.\]  
As in \cite{E} we now introduce formal cohomology 
classes $\delta_{k},W_{k}$ and $\Xi_{k,l}$ as follows. 
Here $\delta_{1},...,\delta_{n}$ 
are formal degree two classes whose sum is zero and whose $r$th elementary  
symmetric polynomial equals $a_{r}$ for $r \geq 2$, and  $W_{k}$ and  
$\Xi_{k,l}$ are given by the formulas 
\begin{equation}\label{wk}W_{k}= 
\sum_{i=1}^{n} f_{i} \frac{\partial \delta_{k}}{\partial 
a_{i}} +\sum_{i=2}^{n} \sum_{j=2}^{n} \sum_{s=1}^{g} 
b_{i}^{s} b_{j}^{s+g}\frac{\partial^{2} 
\delta_{k}}{\partial a_{i} \partial a_{j}},
\end{equation}  
and  
\begin{equation}\label{xikl}\Xi_{k,l}= 
\sum_{s=1}^{g} \left( \sum_{i=2}^{n} b_{i}^{s}\frac{\partial 
\delta_{k}}{\partial a_{i}} -\sum_{j=2}^{n}b_{j}^{s} 
\frac{\partial \delta_{l}}{\partial a_{j}} \right) 
\left(\sum_{i=2}^{n} b_{i}^{s+g} \frac{\partial 
\delta_{k}}{\partial a_{i}}-\sum_{j=2}^{n} b_{j}^{s+g} 
\frac{\partial \delta_{l}}{\partial a_{j}} \right)
\end{equation}  
The tangent bundle $T$ equals $1-\pi_{!}(\mbox{End} V)$ where $V$ is a  
universal bundle over $\mond$ and $\pi: \mond \times\Sigma \to \mond$ 
is the first projection \cite[p.582]{AB}. Arguing as  
in \cite[Prop. 10]{E} we find that the Chern character of $T$ equals  
\begin{equation}\ch T = -\g + \sum_{k=1}^{n} \sum_{l=1}^{n} (\g + W_{l} -W_{k} 
- \Xi_{k,l}) e^{\delta_{k} - \delta_{l}} \label{chT}.\end{equation}  
Let $\lambda_{1}, \ldots ,\lambda_{g}$ denote the  
roots of the equation  
\[x^{g} - x^{g-1} + \frac{1}{2} x^{g-2} - \cdots + \frac{(-1)^{j}}{j!}x^{g-j} 
+ \cdots + \frac{(-1)^{g}}{g!} = 0.\]  
Then \cite[p.862]{K3}  
\[\sum_{j=1}^{g} (\lambda_{j})^{r} = \left\{ \begin{array}{ll} g & r=0\\1 
& r=1\\ 0 & r \geq 2 \end{array} \right. \]  
and hence  
\[\sum_{j=1}^{g} e^{- \lambda_{j} \Xi_{k,l}} = g - \Xi_{k,l}\]  
since $(\Xi_{k,l})^{g+1} = 0$. So expression (\ref{chT}) for 
the Chern character of $T$ may be written as  
\[\ch T =  -\g + \sum_{k=1}^{n} \sum_{l=1}^{n} (W_{l} -W_{k}
 -1) e^{\delta_{k} - \delta_{l}} + \sum_{k=1}^{n} 
\sum_{l=1}^{n}\sum_{j=1}^{g} e^{\delta_{k} - \delta_{l} - \lambda_{j} \Xi_{k,l}},\]  
and using (\ref{newarr}) we find that 
$\sum_{p=0}^{N} t^{p} \ch (\bigwedge{}\!\!^{N-p}\, T^{*})$ equals  
\[(1+t)^{N-n(n-1)g} \prod_{k \neq l} \left( 1 + \frac{e^{\delta_{k} - 
\delta_{l}} - 1}{1+t} \right)^{W_{k}-W_{l}-1} 
\prod_{j=1}^{g} (t +e^{\delta_{l} - \delta_{k} + \lambda_{j} \Xi_{k,l}})\]  
which with some simplifying becomes  
\[(1+t)^{N} \prod_{k \neq l} \left( 1 + \frac{e^{\delta_{k} -\delta_{l}} - 
1}{1+t} \right)^{W_{k}-W_{l}+ \g} \exp 
\left\{\frac{ \Xi_{k,l}}{1 +e^{\delta_{k}-\delta_{l}}t} \right\}.\]  
Now the non-zero Pontryagin roots of $\mond$ are 
$(\delta_{k} - \delta_{l})^{2}$ for $k<l$ and each has 
multiplicity $2g-2$ \cite[Lemma 17]{E}. We also know that  
\[c_{1} = 2 n f_{2} = 2 \sum_{k<l} (W_l - W_k)(\delta_k - \delta_l) -\Xi_{k,l}. \]  
Hence 
\[\td T =  \prod_{k<l} e^{[(W_l - W_k)(\delta_k - \delta_l) -\Xi_{k,l}]} \left( 
\frac{\delta_{k}-\delta_{l}}{2
\mbox{ sinh } \frac{1}{2}(\delta_{k} -\delta_{l})} \right)^{2 g -2}.\]
\indent We have thus shown, with some arranging, that 

\begin{prop}\label{last}
The pairing 
\begin{eqnarray*}\left\{ \prod_{k< l} (t+e^{X_{k,l}})^{\g} 
(t+e^{X_{l,k}})^{\g}\left( \frac{X_{k,l}}{2 
\mbox{\normalshape sinh} \frac{1}{2}(X_{k,l})} 
\right)^{2\g}\left( 1 + \frac{(1-t)(e^{X_{k,l}} - 1)}{1 
+ t e^{X_{k,l}}}\right)^{W_k - W_l} \times \right. \\  
\times \left. \exp \left( \frac{\Xi_{k,l} (1-t^2)}{(1+t 
e^{X_{k,l}})(1+ te^{X_{l,k}})} \right) \right\} [\mond],\end{eqnarray*}  
where $X_{k,l} = \delta_k - \delta_l$, is equal to  
\[\frac{1}{(1+t)^{(n-1)(g-1)}} \left( \prod_{r=1}^{n-1} (
1 - (-t)^{r}) (1-(-t)^{r+1})\right)^{g-1}.\]   
\end{prop}   

\indent Assume now that $n=2$ and let $\alpha = 2f_2,$ $\beta = -4a_2$ and 
$\gamma = 2 \sum_{s=1}^g b_2^s b_2^{s+g}$ denote Newstead's  
generators \cite{N2} for the subring of $H^*({\cal M}_{\Lambda}(2,1))$  
invariant under the induced action of the mapping class group of  
$\Sigma$. In terms of these generators the previous proposition reads  
as (with some simplification),  

\begin{cor}\label{lastt}  When $n=2$ the pairing of  
\[(t+e^X)^{\g} (t+e^{-X})^{\g} \left( \frac{X}{e^{X/2}-e^{-X/2}}
\right)^{2\g} \left( 1 + \frac{(1-t)(e^X-1)}{1 + t e^X} 
\right)^D \exp \left\{ \frac{\xi (1-t^2)}{(1+t e^X)(1+ t e^{-X})} \right\}\]  
with $[\mond]$ equals $(1-t^2)^{g-1}$ where  
\[X = X_{1,2} =  \delta_{1} - \delta_{2} =  \sqrt{\beta}, 
\qquad D =   W_1 - W_2  =  
 -(\alpha \beta + 2 \gamma) 
\beta^{-3/2}, \qquad\xi =  \Xi_{1,2} =   2 \gamma / \beta.\]
\end{cor}

\begin{rem}\indent In the rank two case Thaddeus \cite{T} showed that the  
intersection pairings $\eta [{\cal M}_{\Lambda}(2,1)]$, where 
$\eta\in H^{3g-3}({\cal M}_{\Lambda}(2,1))$, 
can all be expressed in terms of pairings just 
involving $\alpha, \beta$ and $\gamma$. In turn these 
pairings may be expressed in terms of intersection pairings of moduli  
spaces over Riemann surfaces of lower genus involving just $\alpha$  
and $\beta$ by the  equation  
\begin{equation}(\alpha^m \beta^n \gamma^p)[{\cal M}^{g}_{\Lambda}(2,1)] = 
2 g (\alpha^m \beta^n \gamma^{p-1})[{\cal M}^{g-1}_{\Lambda}(2,1)]. 
\label{thad3}\end{equation}.\\
He then determined the pairings 
involving just $\alpha$ and $\beta$ from the Verlinde formula.\\  
\indent Information concerning pairings involving $\alpha$, 
$\beta$and $\gamma$ is contained in 
Corollary \ref{lastt} but is less tractable. 
Not all the pairings involving $\alpha, \beta$ and$\gamma$ may 
be determined from the $\chi(t)$-characteristic but low 
genus calculations suggest that sufficient information is present to  
prove Newstead's and Ramanan's conjecture (that $\beta^g=0$)
\cite[p.344]{N2}. Such a proof would be of 
some interest historically since the$\chi(t)$-characteristic of 
${\cal M}_{\Lambda}(2,1) $ was 
known as early as 1975 \cite{NR}, soon after the Newstead-Ramanan 
conjecture was made, and long before the Verlinde formula was proved 
from which Thaddeus derived his proof.\\  

\end{rem}

\noindent Mathematical Institute, 24-29 St. Giles, Oxford OX1 3LB, England  
\\{\em Email:} earl@@maths.ox.ac.uk, kirwan@@maths.oxford.ac.uk  

\end{document}